\documentclass[]{interact}

\usepackage[utf8]{inputenc}
\usepackage{amsfonts}
\usepackage{amsmath}
\usepackage{amssymb}
\usepackage{mathtools}
\usepackage{amsthm}
\usepackage{color}
\usepackage[final]{microtype}
\usepackage{hyperref}

\usepackage[natbibapa]{apacite}

\usepackage{epstopdf}
\usepackage[caption=false]{subfig}

\usepackage[longnamesfirst,sort]{natbib}
\bibpunct[, ]{(}{)}{;}{a}{,}{,}

\newcommand{\abs}[1]{\left\lvert#1\right\rvert}
\newcommand{\norm}[1]{\left\lvert\left\lvert#1\right\rvert\right\rvert}

\newcommand{\sign}[1]{\mbox{sign}(#1)}

\theoremstyle{plain}
\newtheorem{theorem}{Theorem}[section]
\newtheorem{lemma}[theorem]{Lemma}
\newtheorem{corollary}[theorem]{Corollary}
\newtheorem{proposition}[theorem]{Proposition}
\newtheorem{assumption}[theorem]{Assumption}

\theoremstyle{definition}
\newtheorem{definition}[theorem]{Definition}
\newtheorem{example}[theorem]{Example}

\theoremstyle{remark}
\newtheorem{remark}{Remark}

\begin{document}

\articletype{REGULAR PAPER}

\title{Generating new classes of fixed-time stable systems with predefined upper bound for the settling time\footnote{\textcolor{red}{This is the preprint version of the accepted manuscript: Rodrigo~Aldana-L\'opez, David~G\'omez-Guti\'errez, Esteban Jiménez-Rodríguez, Juan Diego Sánchez-Torres and  Michael Defoort ``Generating new classes of fixed-time stable systems with predefined upper bound for the settling time". International Journal of Control. 2021. DOI: 10.1080/00207179.2021.1936190. 
\textbf{Please cite the publisher's version}. For the publisher's version and full citation details see:
\url{https://doi.org/10.1080/00207179.2021.1936190}.} \textcolor{blue}{
The following links provide access, for a limited time, to a free copy of the publisher's version:
\href{https://www.tandfonline.com/eprint/BMPX5CUSEVFJYHHPMEES/full?target=10.1080/00207179.2021.1936190}{Link 1.} \href{https://www.tandfonline.com/eprint/THCZCUHPMGTTNSPYV7AE/full?target=10.1080/00207179.2021.1936190}{Link 2.}
\href{https://www.tandfonline.com/eprint/XTIW7KNQIJJX6HJGAFTS/full?target=10.1080/00207179.2021.1936190}{Link 3.}
\href{https://www.tandfonline.com/eprint/BQNZBKUT6SAZAGS9VZDT/full?target=10.1080/00207179.2021.1936190}{Link 4.}
}}}

\author{
\name{Rodrigo~Aldana-L\'opez\textsuperscript{a}, David~G\'omez-Guti\'errez\textsuperscript{b,c}, Esteban Jiménez-Rodríguez\textsuperscript{d}, Juan Diego Sánchez-Torres\textsuperscript{e} and  Michael Defoort\textsuperscript{f}
\thanks{CONTACT D.~G\'omez-Guti\'errez. Email: david.gomez.g@ieee.org}}
\affil{\textsuperscript{a}Department of Computer Science and Systems Engineering, University of Zaragoza, Zaragoza, Spain; \textsuperscript{b} Multi-Agent Autonomous Systems Lab, Intel Labs, Intel Tecnología de México, Jalisco, Mexico; \textsuperscript{c} Tecnologico de Monterrey, Escuela de Ingenier\'ia y Ciencias, Jalisco, Mexico; \textsuperscript{d} Relativity6 Inc. Computer Science \& Artificial Intelligence Lab., Jalisco, Mexico; \textsuperscript{e} Research Laboratory on Optimal Design, Devices and Advanced Materials -OPTIMA-, Department of Mathematics and Physics, ITESO, Jalisco, Mexico; \textsuperscript{f} LAMIH, UMR CNRS 8201, INSA, Polytechnic University of Hauts-de-France, Valenciennes, France.}
}

\maketitle

\begin{abstract}
This paper aims to provide a methodology for generating autonomous and non-autonomous systems with a fixed-time stable equilibrium point where an Upper Bound of the Settling Time (\textit{UBST}) is set a priori as a parameter of the system. In addition, some conditions for such an upper bound to be the least one are provided. This construction procedure is a relevant contribution when compared with traditional methodologies for generating fixed-time algorithms satisfying time constraints since current estimates of an \textit{UBST} may be too conservative. The proposed methodology is based on time-scale transformations and Lyapunov analysis. It allows the presentation of a broad class of fixed-time stable systems with predefined \textit{UBST}, placing them under a common framework with existing methods using time-varying gains. To illustrate the effectiveness of our approach, we generate novel, autonomous and non-autonomous, fixed-time stable algorithms with predefined least \textit{UBST}.
\end{abstract}


\section{Introduction}

In recent years, dynamical systems exhibiting convergence to their origin in some finite time, independent of the initial condition of the system, have attracted a great deal of attention. For this class of dynamical systems, their origin is said to be fixed-time stable, which is a stronger notion of finite-time stability~\citep{Bhat2000FiniteTimeSO,Moulay2006}, because in the latter the settling time is, in general, an unbounded function of the initial condition of the system. 
This research effort has derived several contributions on algorithms with the fixed-time convergence property, such as synchronization of complex networks~\citep{Yang2017,Khanzadeh2017,Liu2016finite,Tian2018fixed,Liu2018finite}, stabilizing controllers~\citep{Polyakov2012,Polyakov2015,Basin2016b,Zimenko2018,Sanchez2019IJC,zuo2019fixed,Gomez2020RNC}, distributed resource allocation~\citep{Lin2019}, optimization~\citep{Ning2017distributed}, multi-agent coordination~\citep{Aldana-Lopez2018a,defoort2016fixed,Zuo2014,Wang2018,shi2018new,liu2019finite}, state observers~\citep{Menard2017}, and online differentiation algorithms~\citep{Cruz-Zavala2011,Angulo2013}.

The fixed-time stability property is of great interest in the development of algorithms for scenarios where real-time constraints need to be satisfied. In fault detection, isolation, and recovery schemes~\citep{Tabatabaeipour2014}, failing to recover from the fault on time may lead to an unrecoverable mode. In hybrid dynamical systems, it is frequently required that the observer (resp. controller) stabilizes the observation error (resp. tracking error) before the next switching occurs~\citep{defoort2011,Gomez2015}. In the frequency control of an interconnected power network, not only is the frequency deviation of interest but also how long the frequency stays out of the bounds~\citep{Mishra2018}. 

A Lyapunov differential inequality for an autonomous system to exhibit fixed-time stability was presented in~\citep{Polyakov2012,Zuo2016}, together with an Upper Bound of the Settling Time (\textit{UBST}) of the system trajectory. However, such an upper estimate is too conservative~\citep{Aldana-Lopez2018}. Only recently, non-conservative \textit{UBST} has been derived~\citep{Parsegov2012,Aldana-Lopez2018} for some scenarios. 
An alternative characterization, based on homogeneity theory, was proposed in~\citep{Andrieu2008,Polyakov2016,Tian2018fixedb}. Although it is a powerful tool for the design of high order fixed-time stable algorithms, it poses a challenging design problem for time constrained scenarios, since an \textit{UBST} is often unknown. Thus, the design of fixed-time stable systems where an \textit{UBST} is set a priori explicitly as a parameter of the system, as well as the reduction/elimination of the conservativeness of an \textit{UBST}, is of great interest. This problem has been partially addressed for autonomous systems see, e.g.,~\citep{Aldana-Lopez2018,Sanchez-Torres2018,Aldana-Lopez2018a}, mainly focusing on the class of systems proposed in~\citep{Polyakov2012,Sanchez-Torres2018,Jimenez-Rodriguez2018a}; and for non-autonomous systems, mainly focusing on time-varying gains that either become singular~\citep{Morasso1997,Song2018,Becerra2018,Wang2018,Yucelen2018,Kan2017} or induce Zeno behavior~\citep{Liu2018,Ning2018b} as the predefined-time is reached. 

\textbf{Contributions:} 
We provide a methodology for generating new classes of autonomous and non-autonomous fixed-time stable systems, where an \textit{UBST} is set a priori explicitly as a parameter of the system. The main result is a sufficient condition in the form of a Lyapunov differential inequality, for a nonlinear system to exhibit this property. Additionally, we show that for any fixed-time stable system with continuous settling time function there exists a Lyapunov function satisfying such differential inequality. Based on this characterization, we show how a fixed-time stable system with predefined \textit{UBST} can be constructed from a nonlinear asymptotically stable one, presenting sufficient conditions for such an upper bound to be the least one. 
To illustrate our approach, we provide examples showing how to derive autonomous and non-autonomous fixed-time stable systems, with predefined least \textit{UBST}. This is a significant contribution to the design of control systems satisfying time constraints since, even in the scalar case, the existing \textit{UBST} estimates are often too conservative~\citep{Aldana-Lopez2018}.


\textbf{Notation:} $\mathbb{R}$ is the set of real numbers, $\Bar{\mathbb{R}}=\mathbb{R}\cup\{-\infty,+\infty\}$, $\mathbb{R}_+=\{x\in\mathbb{R}\,:\,x\geq0\}$ and $\Bar{\mathbb{R}}_+=\mathbb{R}_+\cup\{+\infty\}$. The Euclidean norm of $x\in\mathbb{R}^n$ is denoted as $\norm{x}$. $h'(z) = \frac{dh(z)}{dz}$ denotes the first derivative of the function $h:\mathbb{R}\to\mathbb{R}$. 
$C^k(\mathcal{I})$ is the class of functions $f:\mathcal{I}\to\mathbb{R}$ with $k\geq0$ and $\mathcal{I}\subseteq\mathbb{R}$  which has continuous $k$-th derivative in $\mathcal{I}$. For $z \in\mathbb{R}_+$, $\Gamma(z)=\int_0^{+\infty} e^{-\xi}\xi^{z-1}d\xi$ is the Gamma function; for $x,a,b \in\mathbb{R}_+$, $B(x;a,b)=\int_0^x\xi^{a-1}(1-\xi)^{b-1}d\xi$ and $B^{-1}(\cdot;\cdot,\cdot)$ are the incomplete Beta function and its inverse, respectively; for $x \in\mathbb{R}$ $\textnormal{erf}(x)=\int_0^x\frac{2}{\sqrt{\pi}}e^{-\xi^2}d\xi$ is the Error function. 
For $x>0$, $\text{ci}(x) =- \int_x^{+\infty}\frac{\cos(\xi)}{\xi}d\xi$ and $\text{si}(x) =- \int_x^{+\infty}\frac{\sin(\xi)}{\xi}d\xi$ are the cosine and sine integrals respectively.
$\mathcal{K}_a^b$ is the class of strictly increasing $C^1((0,a))$ functions $h:[0,a)\to\Bar{\mathbb{R}}$ with $a,b\in\Bar{\mathbb{R}}$ satisfying  $h(0)=0$ and $\lim_{z\to a}h(z)~=~b$.  

\section{Preliminaries}
\label{Sec:Problem} 

Consider the system
\begin{equation}\label{eq:sys}
    \dot{x}=-\frac{1}{T_c}f(x,t), \ \forall t\geq t_0, \ f(0,t)=0,\ \ 
\end{equation}
where $x\in\mathbb{R}^n$ is the state of the system, $T_c>0$ is a parameter, $t\in[t_0,+\infty)$ and $f:\mathbb{R}^n\times\mathbb{R}_+\to\mathbb{R}^n$, continuous on $x$ and continuous almost everywhere on $t$. The solutions are understood in the sense of Caratheodory~\citep{Regan1997}. We assume that $f(\cdot,\cdot)$
is such that the origin of system~\eqref{eq:sys} is asymptotically stable and system~\eqref{eq:sys} has the properties of existence and uniqueness of solutions in forward-time on the interval $[t_0,+\infty)$~\citep{Khalil2002}. The solution of \eqref{eq:sys} for $t\geq t_0$ with initial condition $x_0$ is denoted by $x(t;x_0,t_0)$, and the initial state is given by $x(t_0;x_0,t_0) = x_0$. 



\begin{remark} For simplicity, throughout the paper, we assume that the origin is the unique equilibrium point of the systems under consideration. Thus, without ambiguity, we refer to the global stability (in the respective sense) of the origin of the system as the stability  of the system. The extension to local stability is straightforward.
\end{remark}

\begin{definition}\citep{Polyakov2014}(Settling-time function)
\label{Def:SettlingTime}
The \textit{settling-time function} of system~\eqref{eq:sys} is defined as
$T(x_0,t_0)=\inf\{\xi\geq t_0:\lim_{t\to\xi}x(t;x_0,t_0)=0\}-t_0\in\Bar{\mathbb{R}}$.
\end{definition}

For autonomous systems ($f$ in~\eqref{eq:sys} does not depend on $t$), the settling-time function is independent of $t_0$. Notice that, if the system is exponentially stable, then according to Definition~\ref{Def:SettlingTime}, $\forall x_0\neq0$, $T(x_0,t_0)=+\infty$.

\begin{definition}\citep{Polyakov2014} \label{def:fixed}(Fixed-time stability) 
System \eqref{eq:sys} is said to be \textit{fixed-time stable} if it is asymptotically stable~\citep{Khalil2002} and the settling-time function $T(x_0,t_0)$ is bounded on  $\mathbb{R}^n\times\mathbb{R}_+$, i.e. there exists $T_{\text{max}}\in\mathbb{R}_+\setminus\{0\}$ such that $T(x_0,t_0)\leq T_{\text{max}}$ if $t_0\in\mathbb{R}_+$ and $x_0\in\mathbb{R}^n$. Thus, $T_{\text{max}}$ is an \textit{UBST} of $x(t;x_0,t_0)$.
\end{definition}

We are interested on finding sufficient conditions on system \eqref{eq:sys} such that an \textit{UBST} is given by the parameter $T_c$, i.e. $T_c=T_{\text{max}}$. 
Of particular interest is to find sufficient conditions such that $T_c$ is the least \textit{UBST}. 



\subsection{Time-scale transformations}

As in~\citep{Pico2013}, the trajectories corresponding to the system solutions are interpreted, in the sense of differential geometry~\citep{Kuhnel2015}, as regular parametrized curves. Since we apply regular parameter transformations over the time variable, then without ambiguity, this reparametrization is sometimes referred to as time-scaling.

\begin{definition}\cite[Definition~2.1]{Kuhnel2015}
\label{Def:RegularParamCurve}
A regular parametrized curve, with parameter $t$, is a $C^1(\mathcal{I})$ immersion $c: \mathcal{I}\to \mathbb{R}$, defined on a real interval $\mathcal{I} \subseteq \mathbb{R}$. This means that $\frac{dc}{dt}\neq 0$ holds everywhere.
\end{definition}

\begin{definition}\cite[Pg.~8]{Kuhnel2015}
\label{Def:RegularCurve}
(Regular parameter transformation)
A regular curve is an equivalence class of regular parametrized curves, where the equivalence relation is given by regular (orientation preserving) parameter transformations $\varphi$, where $\varphi:~\mathcal{I}~\to~\mathcal{I}'$ is $C^1(\mathcal{I})$, bijective and $\frac{d\varphi}{dt}>0$. Therefore, if $c:\mathcal{I}\to\mathbb{R}$ is a regular parametrized curve and $\varphi:\mathcal{I}\to \mathcal{I}'$ is a regular parameter transformation, then $c$  and  $c\circ\varphi:\mathcal{I}'\to\mathbb{R}$ are considered to be equivalent.
\end{definition}


\section{Main Results}
\label{Sec:Main}

The methodology presented in this section to obtain fixed-time stable systems with predefined \textit{UBST} subsumes some existing results, in both, the autonomous and the non-autonomous cases.

\subsection{Time-scaling and settling-time computation}
\label{Subsec:Problem1}

\begin{assumption}
\label{Assump:AsympSys}
$\mathcal{H}:\mathbb{R}\to\mathbb{R}$ is continuous in $\mathbb{R}$, locally Lipschitz in $\mathbb{R}\setminus\{0\}$, and satisfies $\mathcal{H}(0)=0$ and $\forall z\in\mathbb{R}\setminus\{0\}$, $z\mathcal{H}(z)>0$. Moreover, $\tilde{x}(\tau;x_0,0)$ is the unique solution of the asymptotically stable system
\begin{equation}
\label{Eq:x_tau}
\frac{d\tilde{x}}{d\tau}=-\mathcal{H}(\tilde{x}),  \ \  \tilde{x}(0;x_0,0)=x_0, \ \ \tilde{x}\in\mathbb{R},
\end{equation}
and $\mathcal{T}(x_0,0)$ is its settling time.
\end{assumption}

The following lemma presents the construction of the parameter transformation that will be used hereinafter.

\begin{lemma}
\label{Lemma:ParamTransf}
Under Assumption~\ref{Assump:AsympSys}, suppose that
\begin{equation}
\label{Eq:TimeScale}
    \psi(\tau)=T_c\int_{0}^\tau \frac{1}{\Upsilon(\tilde{x}(\xi;x_0,0),\psi(\xi))}d\xi,
\end{equation}
has a unique solution $\psi(\tau)$ on $\mathcal{I}'=[0,\mathcal{T}(x_0,0))$, where $T_c>0$, and $\Upsilon:\mathbb{R}\times\mathbb{R}_+\to\mathbb{R}_+\setminus \{0\}$ is a function such that $[\Upsilon(x,\hat{t})]^{-1}$ continuous for all $x\in\mathbb{R}\setminus \{0\}$ and $\hat{t}\in \mathbb{R}_+$. Then, the map $\psi:\mathcal{I}'\to \mathcal{J}$, where $\mathcal{J}$ is the resulting range of $\psi(\tau)$, is continuous and bijective. Moreover, the bijective function $\varphi:\mathcal{I}=[t_0+\inf\mathcal{J},t_0 + \sup\mathcal{J})\to \mathcal{I}'$  defined by $\varphi^{-1}(\tau) = \psi(\tau) + t_0$  is a parameter transformation. Furthermore, $\mathcal{J}=[0,\lim_{\tau\to\mathcal{T}(x_0,0)}\psi(\tau))$.
\end{lemma}
\begin{proof}
Let $\psi(\tau)$ be the solution of~\eqref{Eq:TimeScale} in $\mathcal{I}'$. Since $\tilde{x}(\tau;x_0,0)$, $\tau\in \mathcal{I}'$, is continuous, then $[\Upsilon(\tilde{x}(\tau;x_0,0),\psi(\tau))]^{-1}$ is continuous on $\tau\in \mathcal{I}'$ and $\psi(\tau)$ is $C^1(\mathcal{I}')$. Moreover, for all $x\in\mathbb{R}\setminus \{0\}$ and $\hat{t}\in \mathbb{R}_+$, $\Upsilon(x,\hat{t})>0$, then it satisfies $\frac{d\psi}{d\tau}>0$, hence $\psi$ is injective \cite[Pg.~34]{Spivak1965}. On the other hand, $\lim_{\tau\to\inf{\mathcal{I}'}}\psi(\tau)=\inf{\mathcal{J}}$ and $\lim_{\tau\to\sup{\mathcal{\mathcal{I}}'}}\psi(\tau)=\sup{\mathcal{J}}$, hence, by the continuity of $\psi$, $\psi$ is surjective and $\mathcal{J} = [0,\lim_{\tau\to\mathcal{T}(x_0,0)}\psi(\tau))$. Thus, $\psi:\mathcal{I}'\to \mathcal{J}$ is bijective. It follows that $\varphi$ is $C^1(\mathcal{I})$, satisfies $\frac{d\varphi}{d\tau}
>0$ and is bijective. Thus, $\varphi:\mathcal{I}\to \mathcal{I}'$ is a parameter transformation. 
\end{proof}

The following lemma shows that if~\eqref{Eq:x_tau} has a known settling time function, then the parameter transformation given in Lemma~\ref{Lemma:ParamTransf} induces a nonlinear system with known settling-time function.

\begin{lemma}\label{Lemma:TimeScale} 
Under Assumption~\ref{Assump:AsympSys}, suppose that $\psi(\tau)$ is the  unique solution of~\eqref{Eq:TimeScale} on $\mathcal{I}'=[0,\mathcal{T}(x_0,0))$, and let
$
\Psi(z,\hat{t}) =\left\lbrace 
\begin{array}{ll}
   \Upsilon(z,\hat{t}) & \text{for } \hat{t}\in\mathcal{J}=[0,\lim_{\tau\to\mathcal{T}(x_0,0)}\psi(\tau)) \\
   1 & \text{otherwise}
\end{array}
\right.
$
with $\Upsilon(z,\hat{t})$ and $T_c$ as in Lemma~\ref{Lemma:ParamTransf}. Then, the system
\begin{equation}
\label{Eq:TSFunc}
    \dot{x}=-\frac{1}{T_c}\Psi(x,\hat{t})\mathcal{H}(x), \ \ x(t_0;x_0,t_0)=x_0,
\end{equation}
where $\hat{t} = t-t_0$, $x\in\mathbb{R}$, has a unique solution 
\begin{equation}
\label{Eq:x_piece}
x(t;x_0,t_0) =\left\lbrace 
\begin{array}{cc}
    \tilde{x}(\varphi(t);x_0,0) & \text{for } t\in \mathcal{I} \\
    0 & \text{elsewhere, } 
\end{array}
\right.
\end{equation}
where $\mathcal{I}=[t_0,t_0 + \lim_{\tau\to\mathcal{T}(x_0,0)}\psi(\tau))$ and $\varphi^{-1}(\tau)=\psi(\tau)+t_0$. Moreover, \eqref{Eq:TSFunc}
is asymptotically stable and the settling time of 
$x(t;x_0,t_0)$
is given by
\begin{equation}
\label{Eq:Time_integral1}    
T(x_0,t_0)=\lim_{\tau\to\mathcal{T}(x_0,0)}\psi(\tau).
\end{equation}
\end{lemma}
\begin{proof}
Consider the parameter transformation $\varphi:\mathcal{I}\to \mathcal{I}'$ given in Lemma~\ref{Lemma:ParamTransf} and let $x(t;x_0,t_0) = \tilde{x}(\varphi(t);x_0,0)$. Notice that, $\tilde{x}(\tau;x_0,0)$, $\tau\in \mathcal{I}'$ and $x(t;x_0,t_0)$, $t\in \mathcal{I}$ are equivalent (in the sense of regular curves). Moreover, $\dot{x} = \frac{d}{dt}\tilde{x}(\varphi(t);x_0,0) = \left.\frac{d\tilde{x}}{d\tau}\right\vert_{\tau=\varphi(t)}\frac{d\varphi}{dt}$, where $\left.\frac{d\tilde{x}}{d\tau}\right\vert_{\tau=\varphi(t)}=-\mathcal{H}(\tilde{x}(\varphi(t);x_0,0))$ and $\frac{d\varphi}{dt} = {T_c}^{-1}\Psi(x(t;x_0,t_0),\hat{t})$. Hence, $\dot{x} = -{T_c}^{-1}\mathcal{H}(x(t;x_0,t_0))\Psi(x(t;x_0,t_0),\hat{t})$, and therefore $x(t;x_0,t_0)$ is solution of \eqref{Eq:TSFunc} on $\mathcal{I}$. Thus, for the solution $\tilde{x}(\tau;x_0,0)$, $\tau\in \mathcal{I}'$ of \eqref{Eq:x_tau}, the equivalent curve $x(t;x_0,t_0)$, $t\in \mathcal{I}$, under the parameter transformation $\varphi$, is solution of \eqref{Eq:TSFunc} on $\mathcal{I}$.

Moreover, since $\frac{d\tilde{x}}{d\tau} = \frac{d}{d\tau}(x\circ\varphi^{-1})(\tau;x_0,0)= -\mathcal{H}(\tilde{x}(\tau;x_0,0))$, then for any solution of \eqref{Eq:TSFunc} on $\mathcal{I}$, there exists an equivalent curve on $\mathcal{I}'$, under the parameter transformation $\varphi^{-1}$, that is solution of \eqref{Eq:x_tau}. Thus the uniqueness of the solution of~\eqref{Eq:TSFunc} on $\mathcal{I}$ follows from the uniqueness of solutions of~\eqref{Eq:x_tau}. 
Finally, since $\tilde{x}(\tau;x_0,0)$ reaches the origin at $\tau=\mathcal{T}(x_0,0)$ then, $x(t;x_0,t_0)$ reaches the origin at $t=t_0 + \lim_{\tau \to\mathcal{T}(x_0,0)}\psi(\tau)$. Moreover, since~\eqref{Eq:TSFunc} has an equilibrium point at $x=0$, then the solution of \eqref{Eq:TSFunc} remains at the origin for all $t\in[t_0 + \lim_{\tau\to\mathcal{T}(x_0,0)}\psi(\tau),+\infty)$. Hence, we can conclude that \eqref{Eq:TSFunc} is asymptotically stable, ~\eqref{Eq:x_piece} is the unique solution in the interval $[t_0,+\infty)$ and the settling time function is given by \eqref{Eq:Time_integral1}.

\end{proof}

\subsection{Fixed-time stability of scalar systems with predefined least \textit{UBST}}

In the rest of the paper, we analyze the cases where $\Psi(x,\hat{t})$ is time invariant or a function only of $t$. We show that in these cases,~\eqref{Eq:TimeScale} has a unique solution.

\begin{assumption}
\label{Assump:Aut}
Let $\Phi:\mathbb{R}_+\to\Bar{\mathbb{R}}_+\setminus\{0\}$ be a function satisfying $\Phi(0)=+\infty$, $\forall z\in\mathbb{R}_+\setminus\{0\}$, $\Phi(z)<+\infty$ and
\begin{equation}
\label{Eq:Finite_Improper}
     \int_0^{+\infty} \Phi(z)dz = 1.
\end{equation}

\end{assumption}

\begin{lemma}
\label{Lemma:ExistenceUniq}
Under Assumption~\ref{Assump:AsympSys}, let $\Upsilon(z,\hat{t})=(\Phi(|z|)\mathcal{H}(|z|))^{-1}$ where  $\Phi(\cdot)$ satisfies Assumption~\ref{Assump:Aut}, then,~\eqref{Eq:TimeScale} has a unique solution on $\mathcal{I}'=[0,\mathcal{T}(x_0,0))$, given by
\begin{equation}
\label{Eq:psiAut}
\psi(\tau)=T_c\int_0^\tau\Phi(|\tilde{x}(\xi;x_0,0)|)\mathcal{H}(|\tilde{x}(\xi;x_0,0)|)d\xi.
\end{equation}

\end{lemma}
\begin{proof}
Let $\Upsilon(z,\hat{t})=(\Phi(|z|)\mathcal{H}(|z|))^{-1}$ and notice that $\Upsilon(\tilde{x}(\tau;x_0,0),\psi)^{-1}$ is independent of $\psi$. Therefore, it follows that
\begin{equation}
\label{Eq:Difpsi}
\frac{d\psi}{d\tau}=T_c\Upsilon(\tilde{x}(\tau;x_0,0),\psi)^{-1},\ \psi(0)=0
\end{equation}
has a unique solution given by $\psi(\tau)=T_c\int_0^\tau\Phi(|\tilde{x}(\xi;x_0,0)|)\mathcal{H}(|\tilde{x}(\xi;x_0,0)|)d\xi$. Moreover, by~\cite[Lemma~1.2.2]{Agarwal1993}, a solution of \eqref{Eq:Difpsi} is also a solution of~\eqref{Eq:TimeScale} and vice versa.
\end{proof}

\begin{assumption}
\label{Assump:NonAut}
Let $\Phi:\mathbb{R}_+\to\Bar{\mathbb{R}}_+\setminus\{0\}$  be a continuous function on $\mathbb{R}_+\setminus\{0\}$ satisfying~\eqref{Eq:Finite_Improper} and $\forall \tau\in\mathbb{R}_+\setminus\{0\}$, $\Phi(\tau)<+\infty$.
Moreover, $\Phi$ is either non-increasing or locally Lipschitz on $\mathbb{R}_+\setminus\{0\}$.
\end{assumption}

\begin{lemma}
\label{Lem:SolNonAut}
Under Assumption~\ref{Assump:AsympSys}, let $\mathcal{I}'=[0,\mathcal{T}(x_0,0))$, and consider the first order ordinary differential equation
\begin{equation}
\label{Eq:NonAutPhi}
\frac{d\psi(\tau)}{d\tau}=T_c\Phi(\tau), \ \ \psi(0)=0,\ \ \tau\in \mathcal{I}',
\end{equation}
where $\Phi(\cdot)$ satisfies Assumption~\ref{Assump:NonAut}, then~\eqref{Eq:NonAutPhi} has a unique solution $\psi:\mathcal{I}'\to \mathcal{J}=[0,\lim_{\tau\to\mathcal{T}(x_0,0)}\psi(\tau))$, which is bijective and given by
\begin{equation}
\label{Eq:NonAutpsi}
\psi(\tau)=T_c\int_{0}^\tau\Phi(\xi)d\xi,  \ \ \tau\in \mathcal{I}'.
\end{equation}
Moreover, let $\Psi(z,\hat{t})=\Phi(\psi^{-1}(\hat{t}))$ then $\psi(\tau)$ is also the unique solution of~\eqref{Eq:TimeScale} on $\mathcal{I}'$.
\end{lemma}
\begin{proof}
It follows that 
\eqref{Eq:NonAutPhi} has a unique solution given by  $\psi(\tau)=T_c\int_{0}^\tau\Phi(\xi)d\xi$. 
Note that $\forall\tau\in \mathcal{I}'\setminus\{0\}$, $\psi:\mathcal{I}'\to \mathcal{J}$ is $C^1(\mathcal{I}')$ with $\frac{d\psi}{d\tau}>0$, and $\psi(0)=0$, hence $\psi$ is injective \cite[Pg.~34]{Spivak1965}. Note that $\lim_{\tau\to\inf{\mathcal{I}'}}\psi(\tau)=\inf{\mathcal{J}}$ and $\lim_{\tau\to\sup{\mathcal{I}'}}=\sup{\mathcal{J}}$ and, by continuity of $\psi$, $\psi$ is surjective. Hence, $\psi:\mathcal{I}'\to \mathcal{J}$ is bijective. Since $\Phi(\psi^{-1}(\psi(\tau)))=\Phi(\tau)$, then $\psi(\tau)$ is a solution of \eqref{Eq:TimeScale}. 
Now, on the one hand if $\Phi$ is non-increasing, then $\Phi\circ\psi^{-1}$ is non-increasing. To show this let $a>b$, $\psi^{-1}(a)>\psi^{-1}(b)$ and $\Phi(\psi^{-1}(a))<\Phi(\psi^{-1}(b))$.
On the other hand, if $\Phi$ is Lipschitz on $[\epsilon,+\infty), \ \forall\epsilon>0$, then $\Phi\circ\psi^{-1}$ is Lipschitz on $J\setminus[0,\epsilon)$. To show this, note that there exists a constant $M_\Phi>0$ such that $|\Phi(\psi^{-1}(x_1))-\Phi(\psi^{-1}(x_2))|\leq M_\Phi|\psi^{-1}(x_1)-\psi^{-1}(x_2)|\leq M|x_1-x_2|$ where $M=M_\Phi\max_{x\in J\setminus[0,\epsilon)}(\psi^{-1})'(x)$.
Then, in the former (resp. in the latter) case it follows from Peano's uniqueness Theorem~\cite[Theorem~1.3.1]{Agarwal1993} (resp. from Lipschitz uniqueness Theorem~\cite[Theorem~1.2.4]{Agarwal1993}) that
\begin{equation}
\label{eq:NonAutEpsilon}
\frac{dz}{d\tau}=T_c\Phi(\psi^{-1}(z)), \ \ z(\epsilon)=T_c\int_0^\epsilon\Phi(\xi)d\xi,
\end{equation}
has a unique solution $z=\psi(\tau)$, $\forall\epsilon,\tau\in \mathcal{I}'\setminus\{0\}$. Since $\psi(0)=0$, then \eqref{eq:NonAutEpsilon} with $\epsilon=0$ has a unique solution $z=\psi(\tau)$, $\tau\in \mathcal{I}'$.
Moreover, by~\cite[Lemma~1.2.2]{Agarwal1993}, a solution of \eqref{eq:NonAutEpsilon} is a solution of~\eqref{Eq:TimeScale} and vice versa.
\end{proof}

In~Lemma~\ref{lemma:characterization_phi}, we present a characterization for a map $\Psi:\mathbb{R}\times\mathbb{R}_+\to\Bar{\mathbb{R}}_+$, in the autonomous case, such that system~\eqref{Eq:TSFunc} is fixed-time stable with $T_c$ as the least \textit{UBST}. 

\begin{lemma}
(Characterization of $\Psi(z,\hat{t})$ for fixed-time stability of autonomous systems with predefined least \textit{UBST})
\label{lemma:characterization_phi}
Under Assumption~\ref{Assump:AsympSys}, with
\begin{equation}
\label{Phi:Aut}
\Upsilon(z,\hat{t}) =(\Phi(|z|)\mathcal{H}(|z|))^{-1} 
\end{equation}
where $\hat{t}=t-t_0$ and $\Phi(\cdot)$ satisfies Assumption~\ref{Assump:Aut}, system~\eqref{Eq:TSFunc} is fixed time stable with $T_c$ as the predefined least \textit{UBST}. 
\end{lemma}

\begin{proof}
By Lemma~\ref{Lemma:ExistenceUniq},~\eqref{Eq:psiAut} is the solution of~\eqref{Eq:TimeScale}. Using the change of variables $z = |\tilde{x}(\tau;x_0,t_0)|$,
\eqref{Eq:Time_integral1} leads to 
$T(x_0,t_0)=T_c\int_0^{\abs{x_0}}\Phi(z)dz
$.
Since $\Phi(\cdot)>0$ then $T(x_0,t_0)$ is increasing with respect to $|x_0|$. Moreover, since $\Phi(\cdot)$ satisfies~\eqref{Eq:Finite_Improper}, the settling-time function satisfies $\sup_{(x_0,t_0) \in \mathbb{R}^n\times\mathbb{R}_+} T(x_0,t_0)=\lim_{\abs{x_0}\to+\infty}T_c\int_0^{\abs{x_0}} \Phi(z)dz=T_c$. 
\end{proof}


The following result states the construction of fixed-time stable non-autonomous systems with predefined \textit{UBST}.

\begin{lemma}\label{Lemma:NAFirst_Order}
(Characterization of $\Psi(z,\hat{t})$ for fixed-time stability of non-autonomous systems with predefined least \textit{UBST})
\label{lemma:NA_first_order}
Under Assumption~\ref{Assump:AsympSys}, let $\psi(\tau)$, $\tau\in \mathcal{I}'=[0,\mathcal{T}(x_0,0))$, be the solution of~\eqref{Eq:NonAutPhi} and $\psi^{-1}(\hat{t})$ its inverse map. Then, with 
\begin{equation}
\label{Phi:NonAut}
\Upsilon(z,\hat{t}) =
    \dfrac{1}{\Phi(\psi^{-1}(\hat{t}))} 
\end{equation}
where $\hat{t}=t-t_0$ and
$\Phi(\cdot)$ satisfies Assumption~\ref{Assump:NonAut}, system~\eqref{Eq:TSFunc} is fixed-time stable with $T_c$ as the predefined \textit{UBST}. Furthermore,
\begin{enumerate}
    \item \label{item1} the settling time is exactly $T_c$ for all $x_0\neq0$ if $\mathcal{T}(x_0,0)~=~+\infty$, for all $x_0\neq 0$;
    \item \label{item2} $T(x_0,t_0)<T_c$ if $\mathcal{T}(x_0,0)<+\infty$, but the least \textit{UBST} is $T_c$ if, in addition, $\mathcal{T}(x_0,0)$ is radially unbounded, i.e. $\mathcal{T}(x_0,0)\to{+\infty}$ as $|x_0|\to{+\infty}$. 
    \item \label{item3} If~\eqref{Eq:x_tau} is fixed-time stable, then, there exists $\Psi_{\mathrm{max}}<+\infty$ such that for all $x_0$ and all $t\in[t_0,t_0+T(x_0,t_0)]$ $\Psi(z,\hat{t})\leq\Psi_{\mathrm{max}}$.
\end{enumerate}
\end{lemma}
\begin{proof}
By Lemma~\ref{Lem:SolNonAut}, the solution of~\eqref{Eq:TimeScale} is given by~\eqref{Eq:NonAutpsi}. Then, the settling time function of \eqref{Eq:TSFunc} is given by $T(x_0,t_0)= T_c\int_0^{\mathcal{T}(x_0,0)}\Phi(\xi)d\xi$. To show item~\textit{(\ref{item1})} note that if $\mathcal{T}(x_0,0)={+\infty}$, then $T(x_0,t_0) = T_c\int_0^{+\infty} \Phi(\xi)d\xi = T_c$, $\forall x_0\in\mathbb{R}\setminus\{0\}$. To show item~\textit{(\ref{item2})} note that, since $\mathcal{T}(x_0,0)<+\infty$ then $T(x_0,t_0)=T_c\int_0^{\mathcal{T}(x_0,0)}\Phi(\xi)d\xi<T_c$. However, if $\mathcal{T}(x_0,0)$ is radially unbounded, then $\sup_{(x_0,t_0) \in \mathbb{R}^n\times\mathbb{R}_+}T(x_0,t_0)=\lim_{|x_0|\to{+\infty}}T_c\int_0^{\mathcal{T}(x_0,0)}\Phi(\xi)d\xi= T_c$. Hence, $T_c$ is the least \textit{UBST}. To show item \textit{(\ref{item3})} notice that, since there exists $\mathcal{T}_{\mathrm{max}}<+\infty$, such that for all $x_0\in\mathbb{R}$,  $\mathcal{T}(x_0,0)\leq\mathcal{T}_{\mathrm{max}}$ then $\exists \hat{T_c}<T_c$, such that $\sup_{(x_0,t_0) \in \mathbb{R}^n\times\mathbb{R}_+}T(x_0,t_0)\leq\lim_{|x_0|\to{+\infty}}T_c\int_0^{\mathcal{T}_{\mathrm{max}}}\Phi(\xi)d\xi=\hat{T}_c$. Thus, for all $t\in[t_0,t_0+\hat{T_c}]$, $\Psi(z,\hat{t})\leq\Psi(z,\hat{T}_c)<+\infty$.
\end{proof}

\begin{remark}
\label{Remark:Zeno}
Fixed-time stability of non-autonomous systems has been applied for the design of stabilizing controllers~\citep{Song2018}, observers~\citep{Holloway2019}, consensus algorithms \citep{Wang2017,Colunga2018b,Wang2018,Ning2018b} and robot control~\citep{Delfin2016} with predefined settling-time at $T_c$, which uses time-varying gains that are either continuous in $[t_0,T_c+t_0)$~\citep{Morasso1997,Song2018,Becerra2018,Wang2018} or piecewise continuous requiring Zeno behavior~\citep{Liu2018,Ning2018b} as $t$ approaches $T_c+t_0$. Notice that, in this paper, we focus on the former case. 
\end{remark}

\begin{remark}
\label{Remark:PsiInfty}
In the autonomous case, $T_c$ is the least \textit{UBST}, whereas, in the non-autonomous case, if item \textit{(1)} is satisfied, every nonzero trajectory converges exactly at $T_c$. This feature has been referred in the literature as predefined-time~\citep{Becerra2018}, appointed-time~\citep{Liu2018} or prescribed-time~\citep{Wang2018}. However, note that $\lim_{t\to t_0+T_c^-}\Psi(z,\hat{t})=+\infty$, but if item \textit{(2)}  or \textit{(3)} 
in Lemma~\ref{lemma:NA_first_order} is satisfied, then the origin is reached before the singularity in $\Psi(z,\hat{t})$ occurs.
\end{remark}

\subsection{Lyapunov analysis for fixed time stability with predefined  \textit{UBST}}
\label{Subsec:Problem2}

The following theorem provides a sufficient condition for a (general) nonlinear system to be fixed-time stable with predefined \textit{UBST}. This result follows from the comparison lemma~\cite[Lemma~3.4]{Khalil2002} and the application of the above results on the time derivative of the Lyapunov candidate function. 

\begin{theorem} (Lyapunov characterization for fixed-time stability with predefined \textit{UBST})
\label{thm:weak_pt} 
Under Assumption~\ref{Assump:AsympSys}, if there exists a continuous and differentiable positive definite radially unbounded function $V:\mathbb{R}^n\to\mathbb{R}$, such that its time-derivative along the trajectories of~\eqref{eq:sys} satisfies
\begin{equation}\label{eq:dV_weak}
\dot{V}(x)\leq-\frac{1}{T_c}\Psi(V(x),\hat{t})\mathcal{H}(V(x)),  \ \  x\in\mathbb{R}^n\setminus\{0\},
\end{equation}  
where $\hat{t} = t-t_0$, and $\Psi(z,\hat{t})$ is characterized by the conditions of Lemma~\ref{lemma:characterization_phi} or  Lemma~\ref{lemma:NA_first_order}, then, system~\eqref{eq:sys} is fixed-time stable with $T_c$ as the predefined \textit{UBST}. If the equality in~\eqref{eq:dV_weak} holds, then $T_c$ is the least \textit{UBST}.
\end{theorem}

\begin{proof}
Let $w(t)$ be a function satisfying $w(t)\geq 0$ and $\dot{w} = -\frac{1}{T_c}\Psi(w,\hat{t})\mathcal{H}(w)$, and let $V(x_0)\leq w(0)$. Then, $T_c$ is the least \textit{UBST} of $w(t)$. Moreover, by the comparison lemma~\cite[Lemma~3.4]{Khalil2002}, it follows that $V(x(t;x_0,t_0))\leq w(t)$. Consequently, $V(x(t;x_0,t_0))$ will converge to the origin before $T_c$. If \eqref{eq:dV_weak} is an equality and $V(x_0)=w(0)$, then, $V(x(t;x_0,t_0)) = w(t)$ and $T_c$ is the least \textit{UBST}. 
\end{proof}

\begin{theorem}
\label{Th:Converse}
If system~\eqref{eq:sys} is autonomous, fixed-time stable and has a continuous settling time function, then there exists a continuous positive definite
function $V:\mathbb{R}^n\to\mathbb{R}$, such that its time-derivative along the trajectories of~\eqref{eq:sys} satisfies
\eqref{eq:dV_weak} with $\Psi(z,\hat{t})$ characterized by the conditions of Lemma~\ref{lemma:characterization_phi}. If in addition, $\lim_{\|x_0\|\to+\infty} T(x_0,t_0)=T_c$
then $V(x)$ is radially unbounded. 
\end{theorem}
\begin{proof}
Let $G(z)=T_c\int_0^z\Phi(\xi)d\xi$, with $\Phi(\cdot)$ satisfying Assumption~\ref{Assump:Aut}. Note that $G'(z)>0,\forall z\geq 0$ and hence $G:\mathbb{R}_+\to[0,T_c)$ is a bijection~\cite[Pg.~34]{Spivak1965}. Moreover note that $G(0)=0$ and $\lim_{z\to\infty}G(z)=T_c$. Hence, $V(x)=G^{-1}(T(x,t_0))$ is a 
continuous and positive definite function satisfying $V(0)=0$. Furthermore, consider the trajectory $x(t;x_0,t_0)$, then, as noted in \cite[Proposition 2.4]{Bhat2000FiniteTimeSO},  $T(x(t;x_0,t_0),t_0)= 
\max\{T(x_0,t_0) - t, 0\}$. 
Therefore, $\dot{V}(x) = -\left(G^{-1}\right)'(T(x,t_0)) =-\frac{1}{T_c}\Phi(V(x))^{-1}= -\frac{1}{T_c}\Psi(V(x),\hat{t})\mathcal{H}(V(x)),\  \forall x\in\mathbb{R}^n\setminus\{0\}$. It follows that, if $\lim_{\|x_0\|\to+\infty} T(x_0,t_0)=T_c$
then $V(x)=G^{-1}(T(x,t_0))$ is radially unbounded. 
\end{proof}


The following theorem allows generating fixed-time stable systems with predefined \textit{UBST} from an asymptotically stable ones that has a Lyapunov function satisfying~\eqref{Eq:LyapInq}. 
By construction, such $V(x)$ will also be a Lyapunov function for system~\eqref{Eq:NonSystemDyn}  satisfying~\eqref{eq:dV_weak}.  

\begin{theorem}
\label{Th.MainResult} (Generating fixed-time stable systems with predefined \textit{UBST})
Under Assumption~\ref{Assump:AsympSys}, let the system
\begin{equation}
\label{Eq:InitSyst}
 \dot{y}=-g(y),
\end{equation}
be asymptotically stable, where $y\in\mathbb{R}^n$, $g:\mathbb{R}^n\to\mathbb{R}^n$ is continuous and locally Lipschitz everywhere except, perhaps, at $y=0$ with $g(0)=0$.
If there exists a Lyapunov function $V(y)$ for system~\eqref{Eq:InitSyst} such that
\begin{equation}
\label{Eq:LyapInq}
    \dot{V}(y)\leq -\mathcal{H}(V(y)),  \ \ \ \forall y\in\mathbb{R}^n,
\end{equation}
then, if $\Psi(V(x),\hat{t})g(x)$ is continuous on $x\in \mathbb{R}^n$,
where $\hat{t} = t-t_0$ and $\Psi(z,\hat{t})$ is a function satisfying the conditions of Lemma~\ref{lemma:characterization_phi} or Lemma~\ref{lemma:NA_first_order}, the system
\begin{equation}
\label{Eq:NonSystemDyn}
\dot{x} = -\dfrac{1}{T_c}\Psi(V(x),\hat{t})g(x), \ \ x(t_0;x_0,t_0)=x_0
\end{equation}
has a unique solution in the interval $[t_0,+\infty)$ and it is fixed-time stable with $T_c$ as the predefined \textit{UBST}.  
\end{theorem}
\begin{proof}
Since the conditions of Lemma~\ref{lemma:characterization_phi} or Lemma~\ref{lemma:NA_first_order} are satisfied, then, \eqref{Eq:TimeScale} has a unique solution. Hence, the proof of the existence of a unique solution for \eqref{Eq:NonSystemDyn} follows by the same arguments as those of the proof of Lemma~\ref{Lemma:TimeScale}.
Let $V(y)$ be a Lyapunov function candidate for~\eqref{Eq:InitSyst} such that \eqref{Eq:LyapInq} holds. Therefore,
$\dot{V}(y) = -\frac{\partial V}{\partial y} g(y)\leq- \mathcal{H}(V(y))$, $\forall y\in\mathbb{R}^n$. Hence, the evolution of $V(x)$ is given by
$\dot{V}(x) = -\frac{1}{T_c}\frac{\partial V}{\partial x}\Psi(V(x),\hat{t})g(x)\leq-\frac{1}{T_c}\Psi(V(x),\hat{t})\mathcal{H}(V(x))$, $\forall x\in\mathbb{R}^n$. Hence, by Theorem~\ref{thm:weak_pt}, $V(x(t;x_0,t_0))$ converges to the origin in fixed-time with $T_c$ as the predefined \textit{UBST}. 
\end{proof}

Notice that,  the term $\Psi(x,\hat{t})\mathcal{H}(x)$ in \eqref{Eq:TSFunc} is continuous at $x=0$ with any choice of $\Psi(x,\hat{t})$ from either Lemma \ref{lemma:characterization_phi} or Lemma \ref{Lemma:NAFirst_Order}, since $\mathcal{H}(0)=0$ and $[\Phi(0)]^{-1}=0$. However, an arbitrary selection of $V(y)$ and $\Psi(z,\hat{t})$ may lead to a right-hand side of~\eqref{Eq:NonSystemDyn} discontinuous at the origin. A construction from a linear system, guaranteeing continuity of the right-hand side of~\eqref{Eq:NonSystemDyn} is provided in the following proposition. 

\begin{corollary}
\label{Prop:ConsLinear}
Let $\Psi(z,t)$ defined as in \eqref{Phi:Aut} with $\Phi(z)$ satisfying Assumption~\ref{Assump:Aut} and $\mathcal{H}(z) = (2\lambda_{\max}(P))^{-1}z$, $P\in\mathbb{R}^{n\times n}$ is the solution of $A^TP + PA=I$ with  $-A\in\mathbb{R}^{n\times n}$ Hurwitz and $\lambda_{\max}(P)$ is the largest eigenvalue of $P$. Then, $\Psi(V(x),\hat{t})Ax$, where $V(x)=\sqrt{x^TPx}$, is continuous, and the system
\begin{equation}
\label{eq:linear_fixed}
\dot{x} = -\frac{1}{T_c}\Psi(V(x),\hat{t})Ax
\end{equation}
where $\hat{t} = t-t_0$, is fixed-time stable with $T_c$ as the predefined \textit{UBST}.
Moreover, if $A=\alpha I+S$ with $\alpha$ a positive constant and $S$ a skew-symmetric matrix then $T_c$ is the least \textit{UBST}.
\end{corollary}
\begin{proof}
Consider system \eqref{Eq:InitSyst} with $g(y)=Ay$ which has a Lyapunov function $V(y)=\sqrt{y^TPy}$ satisfying $\dot{V}\leq- (2\lambda_{\max}(P))^{-1}V(y)=-\mathcal{H}(V(y))$. Note that, $V(\cdot)$ is continuous and $\Psi(\cdot,\hat{t})$ is continuous except at the origin. Therefore, since $\Psi(0,\hat{t})A(0)=0$, to check continuity, it only suffices to show that $\lim_{\|x\|\to 0^+} \|\Psi(V(x),\hat{t})Ax\|=0$ which follows from 
$
\lim_{\|x\|\to 0^+} \|\Psi(V(x),\hat{t})Ax\|^2= 
4\lambda_{\max}(P)^2\lim_{\|x\|\to 0^+}(x^TA^TAx) (x^TPx)^{-1}\Phi(V(x))^{-2}\leq
\frac{4\lambda_{\max}(P)^2\lambda_{\max}(A^TA)}{\lambda_{\min}(P)}\lim_{\|x\|\to 0^+}\Phi(V(x))^{-2} = 0
$. Hence, $\Psi(V(x),\hat{t})Ax$ is continuous everywhere. It follows from Theorem~\ref{Th.MainResult} that \eqref{eq:linear_fixed} is fixed-time stable with $T_c$ as the predefined \textit{UBST}. Note that, if $A = \alpha I+S$ and $P=\frac{1}{2\alpha}I$, then $\dot{V}(y)=-\mathcal{H}(V(y))=-\alpha V(y)$ and $\dot{V}(x) = -\frac{1}{T_c}\Psi(V(x),\hat{t})\mathcal{H}(V(x))$. Hence, by Theorem~\ref{thm:weak_pt}, \eqref{eq:linear_fixed} is fixed-time stable with $T_c$ as the least \textit{UBST}.
\end{proof}


\begin{remark}
Notice that Theorem~\ref{thm:weak_pt} can be used for the design of first and second order controllers as in~\citep{Aldana-Lopez2018}. The design of arbitrary order controllers as in~\citep{Mishra2018}. The design of consensus protocols, as in~\citep{Aldana-Lopez2018a,Ning2017b}; the design of non-autonomous arbitrary order controllers as in~\citep{Pal2020,Gomez2020RNC} or the design of non-autonomous state observers and online differention algorithms~\citep{aldana2020methodology}.
\end{remark}

\section{Examples of fixed-time stable systems with predefined least \textit{UBST}}
\label{Sec:Examples}


\subsection{Examples of autonomous fixed-time stable systems with predefined least \textit{UBST}}

\begin{table*}[t]
    \centering
    \footnotesize{
    \begin{sideways}
    \begin{tabular}{|c|c|c|p{4.5cm}|}
        \hline
        & $\Phi(z)$ & $\dot{x}=-\frac{1}{T_c}(\Phi(\|x\|)\|x\|)^{-1}x$ & Conditions \\
        \hline
        \textit{(i)} & $\frac{1}{\gamma }\left(\alpha h\left(z\right)^p + \beta h\left(z\right)^q\right)^{-k}h'\left(z\right)$ & $\dot{x}= -\frac{\gamma}{T_ch'(\|x\|)}(\alpha h(\|x\|)^{p} + \beta h(\|x\|)^{q})^k\frac{x}{\|x\|}$ & $\gamma =\frac{  \Gamma \left(\frac{1-kp}{q-p}\right) \Gamma \left(\frac{k q-1}{q-p}\right) \left(\frac{\alpha}{\beta}\right)^{\frac{1-kp}{q-p}}}{\alpha^{k}\Gamma (k) (q-p)}$, $kp<1$, $kq>1$, $\alpha,\beta,p,q,k>0$ \\
        \hline
        \textit{(ii)} & $\frac{2}{\pi}\left(\exp\left(2h\left(z\right)\right)-1\right)^{-1/2}h'\left(z\right)$ & $\dot{x}=-\frac{\pi}{2 T_ch'(\|x\|)}(\exp(2h(\|x\|))-1)^{1/2}\frac{x}{\|x\|}$ & $\lim_{z\to 0^+}h'(z)=~+\infty$\\
        \hline
        \textit{(iii)} & $\exp\left(-h\left(z\right)\right)h'\left(z\right)$ & $\dot{x}=-\frac{1}{T_ch'(\|x\|)}\exp(h(\|x\|))\frac{x}{\|x\|}$ &  $\lim_{z\to 0^+}h'(z)=~+\infty$\\
        \hline
        \textit{(iv)} & $\frac{1}{\rho}(\sin(h(z)+\alpha)(1+h(z))^{-2}h'(z)$        & 
        $\dot{x}=-\dfrac{\rho(1+h(\|x\|))^2}{T_ch'(\|x\|)(\sin(h(\|x\|)+\alpha)\|x\|}x$
        &$\rho=\alpha-\text{ci}(1)\cos(1)-\text{si}(1)\sin(1)$, $\alpha>1$ and $\lim_{z\to 0^+}h'(z)=~+\infty$\\
        \hline
    \end{tabular}
    \end{sideways}
    }
    \caption{Examples of $\Phi(z)$ satisfying Assumption~\ref{Assump:Aut}, and fixed-time stable systems with predefined least \textit{UBST} derived from them.}
    \label{Tab:AutonExamp}
\end{table*}

In this subsection, we present the construction of some examples of $\Phi(\cdot)$ satisfying Assumption~\ref{Assump:Aut} for generating autonomous fixed time stable systems with predefined least \textit{UBST}. The result is mainly obtained by applying Corollary~\ref{Prop:ConsLinear}. For simplicity, we take $A=\frac{1}{2}I\in\mathbb{R}^{n\times n}$. 

\begin{proposition}
\label{Proposition:Phis}
Let $h(z)$ be $\mathcal{K}_{+\infty}^{+\infty}$. Then, functions $\Phi(z)$ given in Table~\ref{Tab:AutonExamp}, satisfy Assumption~\ref{Assump:Aut}. Moreover, the system $\dot{x}=-\frac{1}{T_c}(\Phi(\|x\|)\|x\|)^{-1}x$, where $x\in\mathbb{R}^n$ and $-\frac{1}{T_c}(\Phi(\|x\|)\|x\|)^{-1}x$ shown in Table~\ref{Tab:AutonExamp} is fixed-time stable with $T_c$ as the least \textit{UBST}.
\end{proposition}
\begin{proof}
Note that $\frac{1}{\gamma}\int_0^{+\infty} \left(\alpha z^{p} + \beta z^{q}\right)^{-k} dz=\frac{1}{\gamma} \left(\frac{\Gamma \left(\frac{1-kp}{q-p}\right) \Gamma \left(\frac{k q-1}{q-p}\right) \left(\frac{\alpha}{\beta}\right)^{\frac{1-kp}{q-p}}}{\alpha^{k}\Gamma (k) (q-p)}\right)=1$ \citep{Aldana-Lopez2018}, therefore, by Proposition~\ref{col:transformation} then 
$\Phi(z)$ given in Table~\ref{Tab:AutonExamp}-\textit{(i)} satisfies \eqref{Eq:Finite_Improper}. Moreover, since $\Phi(0)=+\infty$, then $\Phi(z)$ satisfies Assumption~\ref{Assump:Aut}. In a similar way, $\Phi(z)$ given in Table~\ref{Tab:AutonExamp}-\textit{(ii)} satisfies Assumption~\ref{Assump:Aut}. The proof that Table~\ref{Tab:AutonExamp}-\textit{(iii)} and Table~\ref{Tab:AutonExamp}-\textit{(iv)} satisfy~\eqref{Eq:Finite_Improper} follows by applying Proposition~\ref{col:transformation} to the functions $F(z)=\exp(-z)$, and $F(z) = \frac{1}{\rho}(\sin(z)+1)(1+z)^{-2}$ (which satisfy $\int_0^{+\infty}F(z)dz$ according to Proposition~\ref{prop:int_beta}), respectively. The proof that the system $\dot{x}=-\frac{1}{T_c}(\Phi(\|x\|)\|x\|)^{-1}x$ with $-\frac{1}{T_c}(\Phi(\|x\|)\|x\|)^{-1}x$ shown in Table~\ref{Tab:AutonExamp} is fixed-time stable with $T_c$ as the least \textit{UBST} follows by applying Corollary~\ref{Prop:ConsLinear} with $g(x)=\frac{1}{2}x$, $V(x)=\|x\|$, $\mathcal{H}(V(x))=\frac{1}{2}\|x\|$ and $\Phi(z)$ given in Table~~\ref{Tab:AutonExamp}.
\end{proof}

\begin{remark}
\label{Remark:Polyakov}
The system in Table~\ref{Tab:AutonExamp}-(i) with $h(z)=z$ reduces to the system analyzed in~\citep{Polyakov2012,Lopez2019}. However, here $T_c$ is given as the predefined least \textit{UBST}. This feature is a significant advantage with respect to~\citep{Polyakov2012,Lopez2019}, because, as illustrated in~\citep{Aldana-Lopez2018}, an \textit{UBST} provided in~\citep{Polyakov2012} is too conservative. Notice that, the fixed-time stable system with predefined \textit{UBST}, analyzed in~\citep{Sanchez-Torres2018}, is found from Table~\ref{Tab:AutonExamp}-(iii) with $h(z) = z^p$ with $0<p<1$. Thus, the algorithms in~\citep{Polyakov2012} and ~\citep{Sanchez-Torres2018} are subsumed in our approach. 
\end{remark}

\begin{figure*}
    \centering
    \def\svgwidth{14.35cm}
\begingroup%
  \makeatletter%
  \providecommand\color[2][]{%
    \errmessage{(Inkscape) Color is used for the text in Inkscape, but the package 'color.sty' is not loaded}%
    \renewcommand\color[2][]{}%
  }%
  \providecommand\transparent[1]{%
    \errmessage{(Inkscape) Transparency is used (non-zero) for the text in Inkscape, but the package 'transparent.sty' is not loaded}%
    \renewcommand\transparent[1]{}%
  }%
  \providecommand\rotatebox[2]{#2}%
  \newcommand*\fsize{\dimexpr\f@size pt\relax}%
  \newcommand*\lineheight[1]{\fontsize{\fsize}{#1\fsize}\selectfont}%
  \ifx\svgwidth\undefined%
    \setlength{\unitlength}{1091.87492311bp}%
    \ifx\svgscale\undefined%
      \relax%
    \else%
      \setlength{\unitlength}{\unitlength * \real{\svgscale}}%
    \fi%
  \else%
    \setlength{\unitlength}{\svgwidth}%
  \fi%
  \global\let\svgwidth\undefined%
  \global\let\svgscale\undefined%
  \makeatother%
  \begin{picture}(1,0.24556383)%
    \lineheight{1}%
    \setlength\tabcolsep{0pt}%
    \put(0,0){\includegraphics[width=\unitlength,page=1]{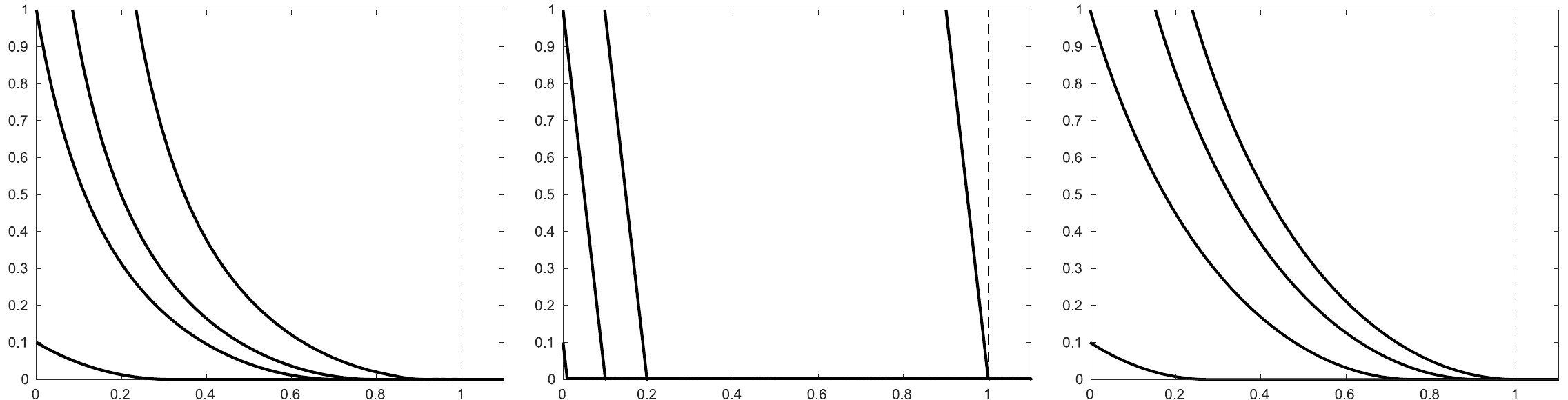}}%
    \tiny{
    \put(0.04946781,0.13316434){\color[rgb]{0,0,0}\rotatebox{-75.672901}{\makebox(0,0)[lt]{\lineheight{1.25}\smash{\begin{tabular}[t]{l}$x_0=0.1$\end{tabular}}}}}%
    \put(0.14505958,0.13314415){\color[rgb]{0,0,0}\rotatebox{-75.546486}{\makebox(0,0)[lt]{\lineheight{1.25}\smash{\begin{tabular}[t]{l}$x_0=1$\end{tabular}}}}}%
    \put(0.17082004,0.13320504){\color[rgb]{0,0,0}\rotatebox{-75.929984}{\makebox(0,0)[lt]{\lineheight{1.25}\smash{\begin{tabular}[t]{l}$x_0=2$\end{tabular}}}}}%
    \put(0.24294369,0.13320505){\color[rgb]{0,0,0}\rotatebox{-75}{\makebox(0,0)[lt]{\lineheight{1.25}\smash{\begin{tabular}[t]{l}$x_0=10e10$\end{tabular}}}}}%
    \put(0.38717776,0.13281919){\color[rgb]{0,0,0}\rotatebox{-73.608826}{\makebox(0,0)[lt]{\lineheight{1.25}\smash{\begin{tabular}[t]{l}$x_0=0.1$\end{tabular}}}}}%
    \put(0.51449294,0.13314416){\color[rgb]{0,0,0}\rotatebox{-75.546486}{\makebox(0,0)[lt]{\lineheight{1.25}\smash{\begin{tabular}[t]{l}$x_0=1$\end{tabular}}}}}%
    \put(0.55570846,0.13320504){\color[rgb]{0,0,0}\rotatebox{-75.929984}{\makebox(0,0)[lt]{\lineheight{1.25}\smash{\begin{tabular}[t]{l}$x_0=2$\end{tabular}}}}}%
    \put(0.57917726,0.13320503){\color[rgb]{0,0,0}\rotatebox{-75}{\makebox(0,0)[lt]{\lineheight{1.25}\smash{\begin{tabular}[t]{l}$x_0=10e10$\end{tabular}}}}}%
    \put(0.71885691,0.13357697){\color[rgb]{0,0,0}\rotatebox{-78.451851}{\makebox(0,0)[lt]{\lineheight{1.25}\smash{\begin{tabular}[t]{l}$x_0=0.1$\end{tabular}}}}}%
    \put(0.7938426,0.13357697){\color[rgb]{0,0,0}\rotatebox{-78.451851}{\makebox(0,0)[lt]{\lineheight{1.25}\smash{\begin{tabular}[t]{l}$x_0=1$\end{tabular}}}}}%
    \put(0.82246312,0.13357697){\color[rgb]{0,0,0}\rotatebox{-78.451851}{\makebox(0,0)[lt]{\lineheight{1.25}\smash{\begin{tabular}[t]{l}$x_0=2$\end{tabular}}}}}%
    \put(0.92710592,0.13379844){\color[rgb]{0,0,0}\rotatebox{-80.142547}{\makebox(0,0)[lt]{\lineheight{1.25}\smash{\begin{tabular}[t]{l}$x_0=10e10$\end{tabular}}}}}%
    \put(0.165,-0.010){\color[rgb]{0,0,0}\rotatebox{0}{\makebox(0,0)[lt]{\lineheight{1.25}\smash{\begin{tabular}[t]{l}time $(t)$\end{tabular}}}}}%
    \put(0.505,-0.010){\color[rgb]{0,0,0}\rotatebox{0}{\makebox(0,0)[lt]{\lineheight{1.25}\smash{\begin{tabular}[t]{l}time $(t)$\end{tabular}}}}}%
    \put(0.84,-0.010){\color[rgb]{0,0,0}\rotatebox{0}{\makebox(0,0)[lt]{\lineheight{1.25}\smash{\begin{tabular}[t]{l}time $(t)$\end{tabular}}}}}%
    }
  \end{picture}%
\endgroup%
    \caption{Examples of autonomous fixed-time systems with $t_0=0$ and $T_c=1$ with $x\in\mathbb{R}$. From left to right: System \eqref{eq:autosys0} with $\alpha = 1,\beta=2,p=0.5,q=2$ and $k=1$; System \eqref{eq:autosys2}; System \eqref{eq:autosys1} with $p=1/2$.}
    \label{fig:SimAutonomous}
\end{figure*}
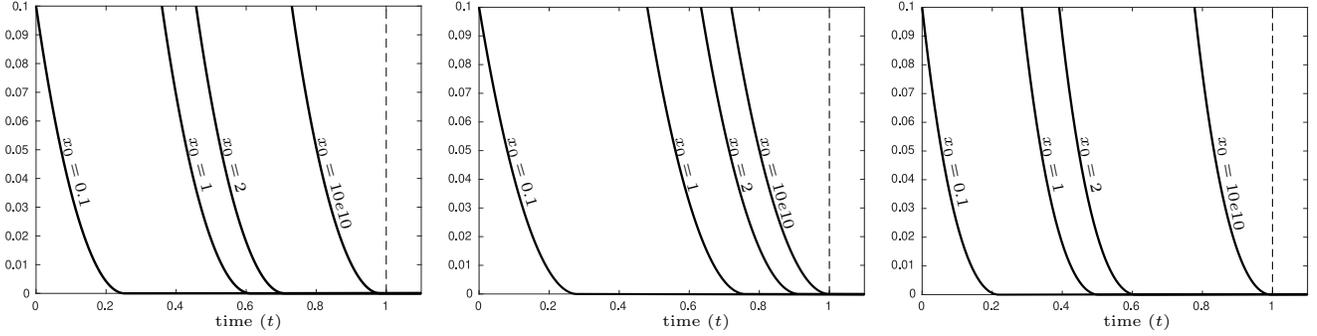

\begin{example}
\label{Exam:Aut}
From Table~\ref{Tab:AutonExamp} new classes of fixed-time stable systems with predefined \textit{UBST}, not present in the literature, can be derived. For instance, the systems
 \begin{equation}
 \label{eq:autosys0}
 \dot{x} = -\frac{\gamma }{T_c}(1+\|x\|)(\alpha\log(1+\|x\|)^p + \beta\log(1+\|x\|)^q)^k\frac{x}{\|x\|}
 \end{equation}
 and
\begin{equation}
 \label{eq:autosys2}
 \dot{x} = -\frac{\pi}{2 T_c}(\exp(2\|x\|)~-~1)^{1/2}\frac{x}{\|x\|}
 \end{equation}
 are obtained from Table~\ref{Tab:AutonExamp}-\textit{(i)} and Table~\ref{Tab:AutonExamp}-\textit{(ii)} with $h(z) = \log(1+z)$ and $h(z)=z$ respectively. 
 Moreover, the system
 \begin{equation}
 \label{eq:autosys1}
 \dot{x} = -\frac{\gamma(\sin(\|x\|^p)+2)}{T_cp(1+\|x\|^p)^2}\frac{x}{\|x\|^p}
 \end{equation}
is obtained from Table~\ref{Tab:AutonExamp}-(iv)  with $h(z) = z^p$ with $0<p<1$. Simulations are shown in Figure~\ref{fig:SimAutonomous}.
\end{example}

\subsection{Examples of non-autonomous fixed-time stable systems with predefined least \textit{UBST}}

In this subsection, we focus on the construction of functions $\Phi(\psi^{-1}(\hat{t}))$ satisfying the conditions of Lemma~\ref{Lemma:NAFirst_Order}. Based on these functions, we provide some examples of non-autonomous systems, with $T_c$ as the settling time for every nonzero trajectory as well as non-autonomous systems with $T_c$ as the least \textit{UBST}.

\begin{table*}
    \centering
    \footnotesize{
    \begin{tabular}{|c|c|p{7.3cm}|}
        \hline
        & $\Phi(\psi^{-1}(\hat{t}))^{-1}$ & Conditions \\
        \hline
          \textit{(i)} &  $T_c\eta'(\hat{t})\abs{1-\eta(\hat{t})}^{-(\alpha+1)}$ & $\alpha\geq0$ and $\eta(z)$ is $\mathcal{K}_{T_c}^1$ and $C^2([0,T_c))$\\
        \hline
        \textit{(ii)} & $\frac{\pi}{2}\sec\left(\frac{\pi\hat{t}}{2T_c}\right)^2\eta'\left(\tan\left(\frac{\pi\hat{t}}{2T_c}\right)\right)$& $\eta(z)$ is $\mathcal{K}_{+\infty}^{+\infty}$ and $C^2([0,+\infty))$\\
        \hline
        \textit{(iii)}& 
        $\frac{\sqrt{\pi}}{2}\eta'\left(\text{erf}^{-1}\left(\frac{\hat{t}}{T_c}\right)\right)\exp\left(\text{erf}^{-1}\left(\frac{\hat{t}}{T_c}\right)^2\right)$ 
        & $\eta(z)$ is $\mathcal{K}_{+\infty}^{+\infty}$ and $C^2([0,+\infty))$\\
        \hline
        \textit{(iv)} & $\gamma(\alpha P(\hat{t})^p + \beta P(\hat{t})^q)^k\eta'(P(\hat{t}))$ & 
        $\gamma =\frac{  \Gamma \left(m_p\right) \Gamma \left(m_q\right) \left(\frac{\alpha}{\beta}\right)^{m_p}}{\alpha^{k}\Gamma (k) (q-p)}$, $kp<1$, $kq>1$, $\alpha,\beta,p,q,k>0$, $m_p=\frac{1-k p}{q-p}$, $m_q=\frac{k q-1}{q-p}$, $\eta(z)$ is $\mathcal{K}_{+\infty}^{+\infty}$, 
        $P(z)^{p-q} = \frac{\beta}{\alpha}B^{-1}\left(\frac{\Gamma(m_p)\Gamma(m_q)z}{\Gamma(k)T_c};m_p,m_q\right)^{-1} -\frac{\beta}{\alpha}$
        \\
        \hline
    \end{tabular}
    }
    \caption{Examples of $\Phi(\psi^{-1}(\hat{t}))$, $\hat{t}=t-t_0$ satisfying conditions of Lemma~\ref{Lemma:NAFirst_Order}, from which non-autonomous fixed-time stable systems with predefined settling time can be constructed. Notice that, in each case, $\lim_{t\to T_c^{-}}\Phi(\psi^{-1}(\hat{t}))^{-1}=+\infty$.}
    \label{Tab:NonAutonExamp}
\end{table*}  
\begin{proposition}\label{lemma:tbg2}
Let $\hat{t} = t-t_0$, then with $\Phi(\psi^{-1}(\hat{t}))$ given in Table~\ref{Tab:NonAutonExamp}, $\Psi(z,\hat{t})$ given in~\eqref{Phi:NonAut}, satisfies the conditions of Lemma~\ref{Lemma:NAFirst_Order}. Moreover, the system  $\dot{x}=-\frac{1}{T_c}\Psi(\|x\|,\hat{t})x$, $x\in\mathbb{R}^n$, is fixed-time stable with $T_c$ as the settling-time for every nonzero trajectory.  
\end{proposition}
\begin{proof}
To show that $\Psi(z,\hat{t})$ given in Table~\ref{Tab:NonAutonExamp}-\textit{(i)} satisfies the condition of Lemma~\ref{Lemma:NAFirst_Order}, choose $h(z) = \frac{1}{\alpha}\left(\frac{1}{(1-\eta(z))^{\alpha}} - 1\right)$. Note that $h(z)$ is $\mathcal{K}_{T_c}^\infty$ for $\alpha\geq 0$. Therefore,  Proposition~\ref{col:tbg_prev} can be used with $\alpha\geq0$. Hence, choosing $\Phi(z)$ as in Proposition~\ref{col:tbg_prev}, leads $ \Phi(\psi^{-1}(\hat{t}))^{-1}=T_ch'(\hat{t}) = \frac{T_c\eta'(\hat{t})}{(1-\eta(\hat{t}))^{\alpha+1}}$. Note that $\psi^{-1}(\hat{t})$ is $C^1([0,T_c))$ and $\eta(\hat{t})$ is $C^2([0,T_c))$, then $\Phi(z)$ is $C^1([0,+\infty))$, therefore satisfies Assumption~\ref{Assump:NonAut}. To show that with $\Phi(\psi^{-1}(\hat{t}))$ given in Table~\ref{Tab:NonAutonExamp}-\textit{(ii)}--Table~\ref{Tab:NonAutonExamp}-\textit{(iv)}, $\Psi(x,\hat{t})$ given in~\eqref{Phi:NonAut}, satisfies the conditions of Lemma~\ref{Lemma:NAFirst_Order}, let $F(z)=\frac{2}{\pi(z^2+1)}$, $F(z)=\frac{2}{\sqrt{\pi}}\exp\left(-z^2\right)$ and $F(z)=\frac{1}{\gamma}(\alpha z^p+\beta z^q)^{-k}$ which satisfies $\int_0^{+\infty}F(z)dz=1$. If $h(z)$ is $\mathcal{K}_\infty^\infty$ and $C^2([0,+\infty))$, by Proposition~\ref{col:transformation}, $\Phi(z)=\frac{2h'(z)}{\pi(h(z)^2+1)}$ and $\Phi(z)=\frac{2h'(z)}{\sqrt{\pi}}\exp\left(-h(z)^2\right)$ satisfy Assumption~\ref{Assump:NonAut}. Furthermore, since  $\Phi(z)=\frac{1}{\gamma}(\alpha h(z)^p+\beta h(z)^q)^{-k}h'(z)$ is non-increasing, it satisfies Assumption~\ref{Assump:NonAut}.  Moreover, by Proposition~\ref{col:transformation},
~\eqref{Eq:NonAutpsi} leads to $\psi(\tau)=\frac{2T_c}{\pi}\arctan(h(\tau))$, $\psi(\tau)=T_c\text{erf}\left(h(\tau)\right)$ and using Proposition \ref{prop:int_beta}, $\psi(\tau)=\frac{T_c(\alpha/\beta)^{m_p}}{\gamma\alpha^k(q-p)} B\left(\left(\frac{\alpha}{\beta}h(\tau)^{p-q}+1\right)^{-1};m_p,m_q\right)$, respectively. Hence, with $\eta(z)=h^{-1}(z)$ and $\eta'(z)=\frac{1}{h'(h^{-1}(z))}$, we obtain $\Phi(\psi^{-1}(\hat{t}))^{-1}$ given in Table~\ref{Tab:NonAutonExamp}-\textit{(ii)}--Table~\ref{Tab:NonAutonExamp}-\textit{(iv)}. From the construction of $\Phi(\psi^{-1}(\hat{t}))^{-1}$ it follows that $\Psi(z,\hat{t})$ satisfies the conditions of Lemma~\ref{Lemma:NAFirst_Order}. The proof that $\dot{x}=-\frac{1}{T_c}\Psi(\|x\|,\hat{t})x$ is a fixed-time stable system with $T_c$ as the settling-time for every nonzero trajectory follows from Lemma~\ref{Lemma:NAFirst_Order}-\textit{(1)}.
\end{proof}

\begin{remark}
\label{Remark:TBGs}
Let $\hat{t} = t-t_0$, then with $\alpha=0$, $\Psi(z,\hat{t})$ in Table~\ref{Tab:NonAutonExamp}-\textit{(i)} reduces to the class of \textit{TBG} proposed in~\citep{Morasso1997}. 
Particular \textit{TBGs}, which can be derived from Table~\ref{Tab:NonAutonExamp}-\textit{(i)}, were used in~\citep{Wang2017,Song2018,Holloway2019,Wang2018,Becerra2018,Colunga2018b,Yucelen2018,Kan2017,Pal2020}. {Notice that, Theorem~1 in\citep{Pal2020} is a particular case of Theorem~\ref{Th.MainResult}, where $\Psi(x,\hat{t})= \frac{1}{1-\hat{t}/T_c}$, $H(z)=\eta(1-e^{-|x|})sign(x)$, and $V(x)=|x|$, with $\hat{t}=t-t_0$, $t_0=0$, and $\eta\geq 1$. Notice that with such $H(z)$, system~\eqref{Eq:x_tau} satisfies, $\mathcal{T}(x_0,0)=+\infty$ for all $x_0\in\mathbb{R}$.
}

\end{remark}

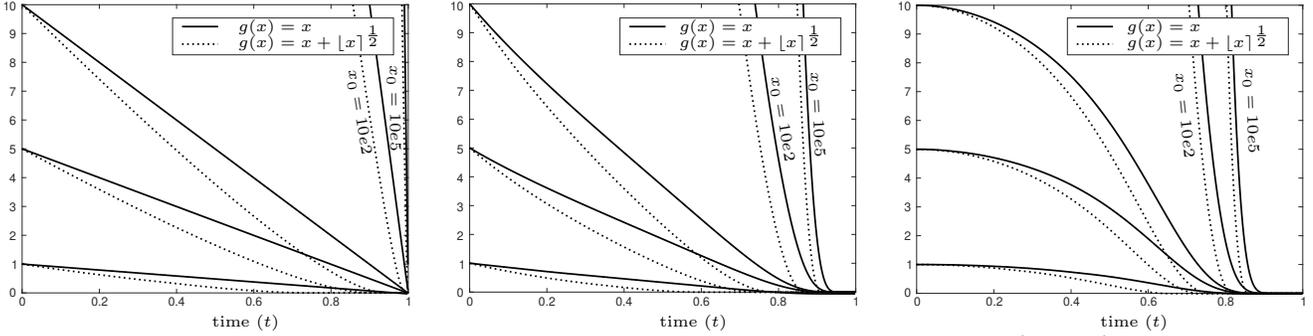
\begin{figure*}
    \centering
    \def\svgwidth{14.35cm}
\begingroup%
  \makeatletter%
  \providecommand\color[2][]{%
    \errmessage{(Inkscape) Color is used for the text in Inkscape, but the package 'color.sty' is not loaded}%
    \renewcommand\color[2][]{}%
  }%
  \providecommand\transparent[1]{%
    \errmessage{(Inkscape) Transparency is used (non-zero) for the text in Inkscape, but the package 'transparent.sty' is not loaded}%
    \renewcommand\transparent[1]{}%
  }%
  \providecommand\rotatebox[2]{#2}%
  \newcommand*\fsize{\dimexpr\f@size pt\relax}%
  \newcommand*\lineheight[1]{\fontsize{\fsize}{#1\fsize}\selectfont}%
  \ifx\svgwidth\undefined%
    \setlength{\unitlength}{1091.87492311bp}%
    \ifx\svgscale\undefined%
      \relax%
    \else%
      \setlength{\unitlength}{\unitlength * \real{\svgscale}}%
    \fi%
  \else%
    \setlength{\unitlength}{\svgwidth}%
  \fi%
  \global\let\svgwidth\undefined%
  \global\let\svgscale\undefined%
  \makeatother%
  \begin{picture}(1,0.24556383)%
    \lineheight{1}%
    \setlength\tabcolsep{0pt}%
    \put(0,0){\includegraphics[width=\unitlength,page=1]{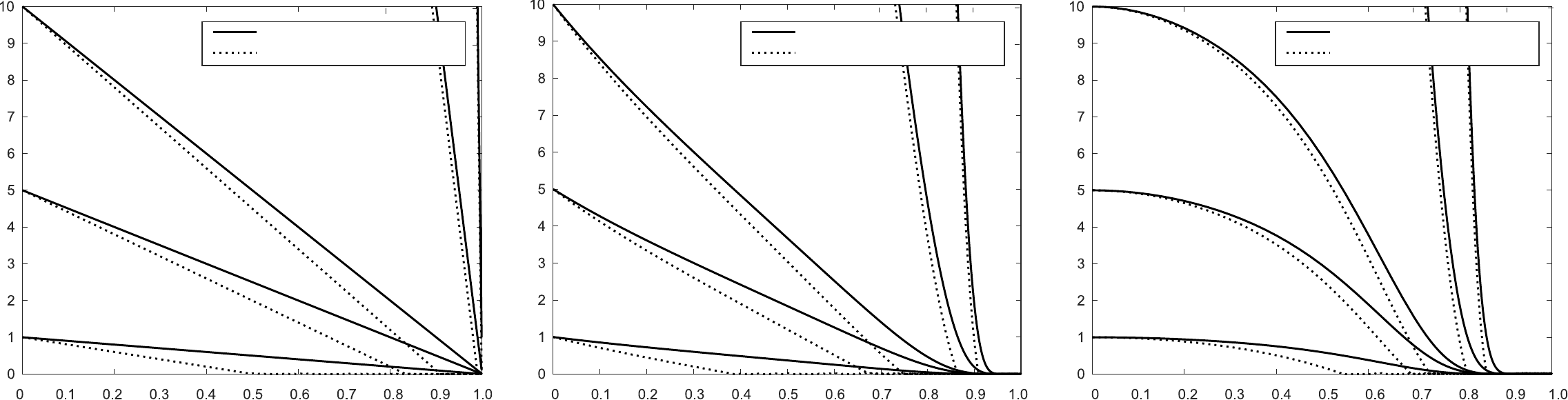}}%
    \tiny{
    \put(0.17,0.21454963){\color[rgb]{0,0,0}\makebox(0,0)[lt]{\lineheight{1.25}\smash{\begin{tabular}[t]{l}$g(x)=x$\end{tabular}}}}%
    \put(0.17,0.2022116){\color[rgb]{0,0,0}\makebox(0,0)[lt]{\lineheight{1.25}\smash{\begin{tabular}[t]{l}$g(x)=x+\lfloor x\rceil^{\frac{1}{2}}$\end{tabular}}}}%
    \put(0.26865448,0.18482757){\color[rgb]{0,0,0}\rotatebox{-80.797663}{\makebox(0,0)[lt]{\lineheight{1.25}\smash{\begin{tabular}[t]{l}$x_0=10e2$\end{tabular}}}}}%
    \put(0.3,0.18874322){\color[rgb]{0,0,0}\rotatebox{-89.356594}{\makebox(0,0)[lt]{\lineheight{1.25}\smash{\begin{tabular}[t]{l}$x_0=10e5$\end{tabular}}}}}%
    \put(0.59,0.17982845){\color[rgb]{0,0,0}\rotatebox{-80.797663}{\makebox(0,0)[lt]{\lineheight{1.25}\smash{\begin{tabular}[t]{l}$x_0=10e2$\end{tabular}}}}}%
    \put(0.62,0.18373915){\color[rgb]{0,0,0}\rotatebox{-86.984048}{\makebox(0,0)[lt]{\lineheight{1.25}\smash{\begin{tabular}[t]{l}$x_0=10e5$\end{tabular}}}}}%
    \put(0.897,0.18292152){\color[rgb]{0,0,0}\rotatebox{-83.810192}{\makebox(0,0)[lt]{\lineheight{1.25}\smash{\begin{tabular}[t]{l}$x_0=10e2$\end{tabular}}}}}%
    \put(0.95,0.18682113){\color[rgb]{0,0,0}\rotatebox{-87.631499}{\makebox(0,0)[lt]{\lineheight{1.25}\smash{\begin{tabular}[t]{l}$x_0=10e5$\end{tabular}}}}}%
    \put(0.51,0.21454963){\color[rgb]{0,0,0}\makebox(0,0)[lt]{\lineheight{1.25}\smash{\begin{tabular}[t]{l}$g(x)=x$\end{tabular}}}}%
    \put(0.51,0.2022116){\color[rgb]{0,0,0}\makebox(0,0)[lt]{\lineheight{1.25}\smash{\begin{tabular}[t]{l}$g(x)=x+\lfloor x\rceil^{\frac{1}{2}}$\end{tabular}}}}%
    \put(0.85,0.21454963){\color[rgb]{0,0,0}\makebox(0,0)[lt]{\lineheight{1.25}\smash{\begin{tabular}[t]{l}$g(x)=x$\end{tabular}}}}%
    \put(0.85,0.2022116){\color[rgb]{0,0,0}\makebox(0,0)[lt]{\lineheight{1.25}\smash{\begin{tabular}[t]{l}$g(x)=x+\lfloor x\rceil^{\frac{1}{2}}$\end{tabular}}}}%
    \put(0.165,-0.010){\color[rgb]{0,0,0}\rotatebox{0}{\makebox(0,0)[lt]{\lineheight{1.25}\smash{\begin{tabular}[t]{l}time $(t)$\end{tabular}}}}}%
    \put(0.505,-0.010){\color[rgb]{0,0,0}\rotatebox{0}{\makebox(0,0)[lt]{\lineheight{1.25}\smash{\begin{tabular}[t]{l}time $(t)$\end{tabular}}}}}%
    \put(0.84,-0.010){\color[rgb]{0,0,0}\rotatebox{0}{\makebox(0,0)[lt]{\lineheight{1.25}\smash{\begin{tabular}[t]{l}time $(t)$\end{tabular}}}}}%
    }
  \end{picture}%
\endgroup%
    \caption{Examples of non-autonomous fixed-time stable system~\eqref{Eq:NonSystemDyn} with $t_0=0$, $T_c=1$, $\lfloor x\rceil^{\frac{1}{2}}=|x|^{\frac{1}{2}}\sign{x}$  and $x\in\mathbb{R}$. From left to right. $\Phi(\psi^{-1}(\hat{t}))^{-1}$ in~\eqref{Syst:ExTBG}; $\Phi(\psi^{-1}(\hat{t}))^{-1}$ in~\eqref{Syst:Sec}; and $\Phi(\psi^{-1}(\hat{t}))^{-1}$ in~\eqref{Syst:poly_eta} with $p=0.5$, $q=2$, $\alpha=1$ and $\beta=2$.}
    \label{fig:SimNonAut}
\end{figure*}

\begin{example}
\label{Example:NonAut}
Let $\hat{t} = t-t_0$, then taking $\eta(z) = z/T_c$ and $\alpha=0$ in~Table~\ref{Tab:NonAutonExamp}-\textit{(i)} results in
\begin{equation}
\label{Syst:ExTBG}
\Phi(\psi^{-1}(\hat{t}))^{-1}= \frac{1}{1-\hat{t}/T_c}
\end{equation}
which corresponds to a \textit{TBG}.
Other time-varying gains, which are not a \textit{TBG} are obtained from Table~\ref{Tab:NonAutonExamp}-\textit{(ii)} and Table~\ref{Tab:NonAutonExamp}-\textit{(iv)} by taking $\eta(z) = z$, which yields to 
\begin{equation}
\label{Syst:Sec}
\Phi(\psi^{-1}(\hat{t}))^{-1} =\frac{\pi}{2}\sec\left(\frac{\pi\hat{t}}{2T_c}\right)^2,
\end{equation}
and
\begin{equation}
\label{Syst:poly_eta}
\Phi(\psi^{-1}(\hat{t}))^{-1} =\gamma(\alpha P(\hat{t})^p + \beta P(\hat{t})^q)^k,
\end{equation}
respectively.
It follows from Lemma~\ref{Lemma:NAFirst_Order}, that taking $g(x)=x$ leads to system~\eqref{Eq:NonSystemDyn} where all non-zero trajectories has $T_c$ as the settling time (Figure~\ref{fig:SimNonAut} solid-line), whereas taking $g(x)=x+\|x\|^{\alpha-1}x$ with $0<\alpha<1$ leads to system~\eqref{Eq:NonSystemDyn}  having $T_c$ as the least \textit{UBST} (Figure~\ref{fig:SimNonAut} dotted-line). Thus, for finite initial conditions, the origin is reached before $T_c$. 
Simulations for the system $\dot{x}=-\frac{1}{T_c}\Psi(|x|,\hat{t})g(x)$ using \eqref{Syst:ExTBG}, \eqref{Syst:Sec} and \eqref{Syst:poly_eta}, $x\in\mathbb{R}$ are presented in Figure~\ref{fig:SimNonAut}.
\end{example}


\subsection{Examples of fixed-time second order systems with predefined \textit{UBST}}
\label{Sec:Exam2ndOrder}

\begin{proposition}
\label{Eq:DiffError}
Assume that, under a suitable selection of $k_1$, $k_2$, $g_1(\cdot)$ and $g_2(\cdot)$, the system
\begin{align}
  \label{Eq:Diff}
  \dot{y}_1&=-k_1g_1(y_1)+y_2\\
  \dot{y}_2&=-k_2g_2(y_1)-y_2
\end{align}
is finite-time stable and there exists a Lyapunov function $V(y)$, satisfying \eqref{Eq:LyapInq}.
Then, the system
\begin{align}
\label{Eq:PredefDiff}
  \dot{z}_1&=-\kappa(\hat{t})k_1g_1(z_1)+z_2\\
  \dot{z}_2&=-\kappa(\hat{t})^2k_2g_2(z_1)
\end{align}
where $\kappa(\hat{t})=\frac{1}{T_c}\Psi(x,\hat{t}) $ with $\Psi(x,\hat{t})$ given as in~\eqref{Phi:NonAut} with $\Phi(\psi^{-1}(\hat{t}))^{-1}$ given in~\eqref{Syst:ExTBG} and $\hat{t} = t-t_0$, is fixed-time stable with $T_c$ as the predefined \textit{UBST}.
\end{proposition}
\begin{proof}
Consider the coordinate change $z_1=x_1$ and $z_2=\kappa(\hat{t})x_2$. Then, in the new coordinates the dynamic is represented by
$\dot{x}_1=\frac{1}{T_c}\Psi(x,\hat{t})[-k_1g_1(x_1)+x_2]$ and $\dot{x}_2=\frac{1}{T_c}\Psi(x,\hat{t})[-k_2g_2(x_1)-x_2]$.
Hence, the result follows from Theorem~\ref{Th.MainResult} by 
taking
$g(y)=
\left[
     k_1g_1(y_1)-y_2, \
     k_2g_2(y_1)+y_2
\right]^T$.
\end{proof}

\begin{proposition}
\label{Prop:Control}
Assume that, under a suitable selection of $k_1$, $k_2$, $g_1(\cdot)$ and $g_2(\cdot)$, the system
\begin{align}
  \label{Eq:Control}
  \dot{y}_1&=y_2\\
  \dot{y}_2&=-k_1g_1(y_1)-k_2g_2(y_2)-y_2
\end{align}
is finite-time stable and there exists a Lyapunov function $V(y)$, satisfying \eqref{Eq:LyapInq}.
Then, the system 
\begin{align}
    \dot{z}_1&=z_2\\
    \dot{z}_2&=-\kappa(\hat{t})^2k_1g_1(z_1)-k_2\kappa(\hat{t})^2g_2(\kappa(\hat{t})^{-1}
 z_2)
 \label{Eq:PredControl}
\end{align}
where $\kappa(\hat{t})=\frac{1}{T_c}\Psi(x,\hat{t}) $ with $\Psi(x,\hat{t})$ given as in~\eqref{Phi:NonAut} with $\Phi(\psi^{-1}(\hat{t}))^{-1}$ given in~\eqref{Syst:ExTBG} and $\hat{t} = t-t_0$, is fixed-time stable with $T_c$ as the predefined \textit{UBST}. 
\end{proposition}
\begin{proof}
The proof is similar to the one given for Proposition~\ref{Eq:DiffError}, considering the coordinate change $z_1=x_1$ and $z_2=\kappa(\hat{t})x_2$. 

\end{proof}

\begin{remark}
Notice that the result in Proposition~\ref{Eq:DiffError} can be applied straightforwardly to the design of predefined-time second-order observers, whereas the result in Proposition~\ref{Prop:Control} can be applied straightforwardly to the design of second-order predefined-time controllers. These results can be extended to the high order case.
\end{remark}

\begin{remark}
An important consequence of Lemma~\ref{Lemma:NAFirst_Order} is that, based on the settling-time function of~\eqref{Eq:PredefDiff} (resp.~\eqref{Eq:PredControl}), $T_c$ can be the least \textit{UBST}. Moreover, If~\eqref{Eq:Diff} (resp.~\eqref{Eq:Control}) is fixed-time stable, then $\kappa(\hat{t})$ is bounded for all $t\in[t_0,t_0+T(x_0,t_0)]$ and all $x_0\in\mathbb{R}^2$.
\end{remark}

\section{Conclusions and future work}
\label{Sec:Conclusions}
We presented a methodology for generating fixed-time stable algorithms such that an \textit{UBST} is set a priori explicitly as a parameter of the system, proving conditions under which such upper bound is the least one. Our analysis is based on time-scaling and Lyapunov analysis. We have shown that this approach subsumes some existing methodologies for generating autonomous and non-autonomous fixed-time stable systems with predefined \textit{UBST} and allows to generate new systems with novel vector fields. Several examples are given showing the effectiveness of the proposed method. As future work, we consider the application/extension of these results to differentiators, control and consensus algorithms.

\appendix
\section{Auxiliary identities}
\begin{proposition}\label{prop:int_beta}
The following identities are satisfied: 
$
i)\int_0^{+\infty}\frac{\sin(z)+a}{(1+z)^2}dz = a-\rm{ci}(1)\cos(1)-\rm{si}(1)\sin(1)$; 
$
 ii)\int_0^{x}(\alpha z^p + \beta z^q)^{-k}dz=\frac{(\alpha/\beta)^{m_p}}{\alpha^k(q-p)} B\left(\left(\frac{\alpha}{\beta}x^{p-q}+1\right)^{-1};m_p,m_q\right)$,
for $kp<1$, $kq>1$, $\alpha,\beta,p,q,k>0$, $m_p=\frac{1-k p}{q-p}$, $m_q=\frac{k q-1}{q-p}$ and $a>1$.
\end{proposition}
\begin{proof}
\textit{i)} It follows from $\int_0^{+\infty}\frac{a}{(1+z)^2}dz = a$ and the change of variables $u=1+z$ with integration by parts and the definition of $\text{ci}(z)$ and $\text{si}(z)$. \textit{ii)} It follows by the change of variables $u=\left(\frac{\alpha}{\beta}x^{p-q}+1\right)^{-1}$ using the definition of $B(\cdot;\cdot,\cdot)$ similarly as in \citep{Aldana-Lopez2018}. 
\end{proof}


\section{Some results on the construction of \lowercase{$\Phi(z)$}}
\label{Append:FindingPhi}




\begin{proposition}\label{col:transformation}
Let $h(\cdot)$ be a $\mathcal{K}_{+\infty}^{+\infty}$ function and let $F:\mathbb{R}_+\to\Bar{\mathbb{R}}\{0\}$ a function satisfying $\int_0^{+\infty}F(z)dz=M$. Then, $\Phi(z) = \frac{1}{M}F(h(z))h'(z)$ satisfies \eqref{Eq:Finite_Improper}.
Furthermore, with such $\Phi(z)$, \eqref{Eq:NonAutpsi} becomes $\psi(\tau)=T_c\int_0^{h(\tau)}F(\xi)d\xi$. Moreover, if $F(z)<+\infty, \forall z\in\mathbb{R}_+$ and $\lim_{z\to 0^+}h'(z)=+\infty$, then $\Phi(\cdot)$ satisfies Assumption~\ref{Assump:Aut}. 
\end{proposition}
\begin{proof}
Using $\xi = h(z)$, it follows
$\int_0^{+\infty} \Phi(z)dz = \frac{1}{M}\int_0^{+\infty} F(h(z))h'(z)dz = \frac{1}{M}\int_0^{+\infty} F(\xi)d\xi = 1$. Moreover, if $F(z)<+\infty, \forall z\in\mathbb{R}_+$ and $\lim_{z\to 0^+}h'(z)=+\infty$, then $\lim_{z\to 0^+}\Phi(z) = F(0)\lim_{z\to 0^+}h'(z)=+\infty$. Hence, $\Phi(\cdot)$ satisfies Assumption~\ref{Assump:Aut}. The rest of the proof follows from \eqref{Eq:NonAutpsi} and the change of variables $u=h(\xi)$.
\end{proof}

\begin{proposition}\label{col:tbg_prev}
Let $h(z)$ be a $\mathcal{K}_{T_c}^{\infty}$ function. Then, the function $\Phi(z)$ characterized by
$\Phi(h(z)) = \frac{1}{T_c}\left(\frac{dh(z)}{dz}\right)^{-1}$
satisfies Assumption~\ref{Assump:NonAut}. Moreover, let $\Upsilon(z,\hat{t})=\Phi(\psi^{-1}(\hat{t}))$ then, the solution of \eqref{Eq:TimeScale} is given by $\psi(\tau)=h^{-1}(\tau)$ and $\Phi(\psi^{-1}(\hat{t})) = \Phi(h(\hat{t}))$.
\end{proposition}
\begin{proof}
Using the change of variables $\xi=h(z)$,  $\int_0^{+\infty}\Phi(\xi)d\xi=\int_0^{T_c}\Phi(h(z))\frac{dh}{dz}dz = \frac{1}{T_c}\int_0^{T_c}dz=1$.
Then, from Lemma~\ref{Lem:SolNonAut}, \eqref{Eq:NonAutpsi} is the solution of \eqref{Eq:TimeScale}. Moreover, $\psi(\tau) = T_c\int_0^{h^{-1}(\tau)}\Phi(h(z))\frac{dh}{dz}dz= h^{-1}(\tau)$. Hence, $\psi^{-1}(\hat{t})=h(\hat{t})$ and $\Phi(\psi^{-1}(\hat{t})) = \Phi(h(\hat{t}))$.
\end{proof}





\end{document}